\input amstex
\documentstyle{amsppt}
\magnification=\magstep1
%\NoBlackBoxes
\TagsOnRight

\document

\topmatter
\title
Periodic points
of quasianalytic Hamiltonian billiards
\endtitle
\author
M\. Novitskii and Yu\. Safarov
\endauthor
\date
January 1999
\enddate
\address
Institute for Low Temperature Physics, Kharkov
\endaddress
\email
novitskii\@ilt.kharkov.ua
\endemail
\address
{\noindent Department of Mathematics, King's College,
Strand, London WC2R 2LS, UK}
\endaddress
\email
ysafarov\@mth.kcl.ac.uk
\endemail
\thanks
The first author would like to thank the Department
of Mathematics of King's College, London, where this work
was performed, for its  hospitality.
The first author was supported by the Royal Society and the
second author by the British Engineering and Physical
Sciences Research Council, Grant B/93/AF/1559. Both
authors are very grateful to D\. Vassiliev for valuable
discussions.
\endthanks
\abstract
We study absolutely periodic points and trajectories
of Hamiltonian systems. Our main result is
a necessary and sufficient  for a Hamiltonian system
to have the following property: if there exists
one absolutely periodic trajectory then all
trajectories are periodic with the same period.
\endabstract
\endtopmatter

\def\dr{\text{\rm d}}

\head
{Introduction}
\endhead

Periodic and absolutely periodic points of
Hamiltonian systems play an important role
in the study of spectral and semi-classical
asymptotics for partial differential operators
(see, for example, \cite{DG, PP2, SV2}).
For instance, if the set of periodic points
of the geodesic flow on a compact Riemannian manifold
is of measure zero then the classical Weyl two-term
asymptotic formula holds, which means, roughly speaking,
that the eigenvalues of $\sqrt{-\Delta}$ ($\Delta$ being
the Laplace-Beltrami operator) are distributed uniformly
on the real line. If this set has a positive
measure then the spectrum of $\sqrt{-\Delta}$ may contain
clusters --- contracting
groups of eigenvalues of unusually high total
multiplicity \cite{DG, SV2}.
Similar results have been also obtained for
\roster
\item"(i)"
general
self-adjoint elliptic partial differential operators
on manifolds without boundary,
where one has to consider the Hamiltonian flow
generated by the principal symbol \cite{DG, SV2};
\item"(ii)"
boundary value
problems, where instead of the Hamiltonian flows
one has to deal with so-called Hamiltonian billiards
\cite{I, SV2};
\item"(iii)"
semi-classical problems, where the principal
symbol is defined in a different way and is not
necessarily a homogeneous function \cite{PP2}.
\endroster

General two-terms asymptotic formulae for
the counting function of a partial differential
operator contain an oscillating term, given
by an integral over the set of periodic points
and involving the period function
\cite{PP2, SV2}. These formulae often imply
that the spectrum does contain
clusters whenever the period function is constant
\cite{SV2}. Therefore it is essential
to know
\roster
\item"(i)" whether the set of
periodic points is of measure zero and, if not,
\item"(ii)" whether the period function is
constant.
\endroster
Unfortunately, there are very few results on these
problems. Even in the simplest case of the standard
Euclidean billiard in a convex smooth domain it is
still unknown whether the set of periodic points
is of measure zero.

The main aim of this paper is to show how one
can obtain results on the above mentioned problems
with the use of theory of quasianalytic functions.
We shall take advantage of the notion of
absolutely periodic point introduced in \cite{DG}.
An absolutely periodic point can be
roughly described as a periodic point at which all
the derivatives of the Hamiltonian flow coincide with
that of the identical map (see Definition 4.1).

The set of periodic points may have a very
complicated structure. In the generic case,
the set of absolutely periodic points is much
poorer. However, it is known that these two sets are of
the same measure (see Lemma 4.4). Therefore it is tempting
to find conditions under which there
are no absolutely periodic points and, consequently, the
set of periodic points is of measure zero.

Sometimes such conditions can be formulated in terms of
smoothness. In \cite{SV2} the authors considered, in
particular, the Hamiltonian systems generated by
analytic homogeneous Hamiltonians.
It was proved that in the analytic case the existence of one
absolutely periodic point implies that all sufficiently
close points are absolutely periodic. This allows one
to find sufficient conditions for the set of periodic
points to be empty.

In this paper we use the scale of Carleman spaces, which
contains the class of analytic functions as a particular
case. This scale can be divided into two parts:
quasianalytic and nonquasianalytic classes. Our main
observation is that the technique suggested in
\cite{SV2} is also applicable in the quasianalytic case.
We extend results obtained in \cite{SV2} to the
quasianalytic Carleman classes
and construct counterexamples which imply that
these results do not hold for nonquasianalytic
Hamiltonians. Thus, we obtain necessary and
sufficient conditions in terms of the Carleman spaces.

The paper is organized as follows.
In Section 1 we define the Carleman classes
and recall some results, including
a criterion of quasianalyticity.
Section 2 and 3 are devoted to definitions and
basic results concerning Hamiltonian flows
and billiards. In Section 4 we define
periodic and absolutely periodic points and
trajectories and formulate our main results
(note that Lemma 4.2 is new even in the analytic case).
These results are proved in Section 6. In Section
5 we give an explicit description of geodesics on a
surface of revolution, which we use later
in Section 5 in order to construct the counterexamples.

As we have mentioned before, the results on
periodic and absolutely points have applications
in spectral theory. Without getting
into details, we formulate a simple corollary of
our Corollary 4.6 and \cite{SV2, Theorem~1.6.1}.

\proclaim{Corollary}
Let $M$ be a compact quasianalytic Riemannian
manifold with strictly convex boundary and $N(\lambda)$ be
the counting function of $\sqrt{-\Delta}\,$, where $\Delta$
is the Laplace-Beltrami operator subject to Dirichlet or
Neumann boundary condition. Then
$$
N(\lambda)=c_0\lambda^d+c_1\lambda^{d-1}+o(\lambda^{d-1})\,,
\qquad\lambda\to+\infty\,,\tag0.1
$$
where $c_0$ and $c_1$ are the classical Weyl
coefficients. In particular, {\rm (0.1)} is valid for
the usual Laplacian on a convex domain with
quasianalytic strictly convex boundary.
\endproclaim

\head
{1. Carleman spaces $ C(m_n)$}
\endhead

Let $\{m_n\}_{n=1}^\infty$ be a nondecreasing sequence
of positive numbers satisfying the following
conditions:
\roster
\item
$m_n$ grow faster than any power of $n$;
\item
$m_{n+1}/m_{n}\le n\,C^n$ for some positive constant $C\,$;
\item
the sequence $\{\ln(m_n/n!)\}_{n=1}^\infty$ is convex
with respect to $n\,$.
\endroster

\definition{Definition 1.1}
Let $D\subset\Bbb R^d$ be a bounded domain.
The Carleman space $C(m_{n},D)$ is defined as the
set of all functions $u\in C^\infty(D)$ such that
$$
\sup_{x\in D}|\partial^\alpha u|
\le C^{|\alpha|}\,m_{|\alpha|}\,,\qquad\forall\alpha\,,
$$
where $\partial^\alpha u$ are the partial derivatives of $u$
and $C$ is a positive constant which may depend on $u$
but is independent of $\alpha$. We say that $u\in C(m_n)$
if $u\in C(m_{n},D)$ for any bounded open subset $D$ of its
domain of definition.
\enddefinition

\example{Example 1.2}
If $\,m_n=n!\,$ then $C(m_n,D)$ coincides with the class
of real analytic functions.
\endexample

\example{Example 1.3}
If $s\ge1$ then the sequence $m_n = n^{sn}$ satisfy the
conditions (1)--(3). The space $C(m_n)$ defined by such
a sequence is called the Gevrey class.
\endexample

In a similar way one defines classes $C(m_n)$
of vector-valued and matrix-valued functions.
We shall take advantage of the following results
(see \cite{D, Section~5}).

\proclaim{Lemma 1.4}
Composition of $C(m_n)$-functions is a $C(m_n)$-function.
\endproclaim

\proclaim{Lemma 1.5}
Let $F(t,z)$ be a real $C(m_n)$-function of variables
$t\in\Bbb R\,$, $z\in\Bbb R^N$. If
$\partial_t F(t_0,z_0)\ne0$ then the local $t$-solution
$t^*(z)$ of the equation $F(t,z)=0$, defined in a
neighbourhood of $z_0$, belongs to $C(m_n)$.
\endproclaim

\proclaim{Lemma 1.6}
If the Cauchy problem for the
a system of ordinary differential equations
$$
\frac{\partial z}{\partial t}=f(t,z(t))\,,
\qquad z(0)=z_0\in\Bbb R^N\,,
$$
has a continuously differentiable solution
$z(t;z_0)$ and $f\in C(m_n)$ then
$z(t;z_0)$ also belongs to $C(m_n)$ as a
function of $(t,z_0)$.
\endproclaim

In view of Lemma 1.4 one can define $C(m_n)$-manifolds and
$C(m_n)$-functions on $C(m_n)$-manifolds.
The scale of spaces $C(m_n)$ on real analytic manifolds
was used for the study of the wave front sets of solution
of linear differential  equations  with
$C(m_{n})$-coefficients \cite{H, Chapter~8} and
some problems of function theory \cite{D}.

\remark{Remark 1.7}
One can define the classes $C(m_n)$ assuming only that
$m_n$ grow faster than any power of $n\,$. The sequences
$\{m_n\}$ satisfying the conditions (2) and (3) are
said to be regular. Lemmas 1.4--1.6 were proved in \cite{D}
only for the regular sequences. It is quite possible
that these lemmas hold true under less restrictive  conditions
(for instance, one can try to apply the regularization theory
\cite{M}),
and then all our further results remain valid as well.
\endremark

\definition{Definition 1.8}
The class $C(m_n)$ is said to be {\it quasianalytic} if
every function $u\in C(m_n)$ which has an infinite order
zero is identically equal to zero.
\enddefinition

\example{Example 1.9}
The analytic functions are quasianalytic.
\endexample

The following theorem gives a criterion of quasianalyticity
for the classes $C(m_n)$.

\proclaim{Theorem 1.10 \cite{M}}
$C(m_n)$ is quasianalytic if and only if the
series $\sum_{n=1}^{\infty}m_n^{-1/n}$ is divergent.
\endproclaim

\head
{2. Hamiltonian flows}
\endhead

\subhead
{1. Hamiltonian flows in $\Bbb R^d\times\Bbb R^d$}
\endsubhead
Let $h(x,\xi)$ be a real infinitely differentiable function on
$\Bbb R^d\times\Bbb R^d$, where $x=(x_1,\ldots,x_d)$ and
$\xi=(\xi_1,\ldots,\xi_d)$ are $d$-dimensional
variables. The systems of ordinary differential equations
$$
\alignat2
\dot x^*(t;y,\eta)
&=h_\xi\bigl(x^*(t;y,\eta),\xi^*(t;y,\eta)\bigr),\qquad
&&x^*(0;y,\eta)=y\,,\\
\dot\xi^*(t;y,\eta)&
=-h_x\bigl(x^*(t;y,\eta),\xi^*(t;y,\eta)\bigr),\qquad
&&\xi^*(0;y,\eta)=\eta\,,
\endalignat
$$
is said to be a {\sl Hamiltonian system\/}, and the function
$h$ is called the {\sl Hamiltonian\/}. If the first and
second derivatives of $h$ are uniformly bounded then
the Hamiltonian system has a unique global
solution $(x^*,\xi^*)$ for every set of initial
data $(y,\eta)\in\Bbb R^d\times\Bbb R^d$ (see, for example,
\cite{Ha}).

A solution
$$
\bigl(x^*(t),\xi^*(t)\bigr)\ :=\
\bigr(x^*(t;y,\eta),\xi^*(t;y,\eta)\bigl)\,,
\qquad t\in\Bbb R,\tag2.1
$$
of the Hamiltonian system
is usually interpreted as a trajectory in the {\sl phase space}
$\Bbb R^{2d}=\Bbb R^d\times\Bbb R^d$ originating from the initial point
$(y,\eta)$. The group of shifts along the Hamiltonian
trajectories is said to be the {\it Hamiltonian flow\/}.
It is easy to see that the Hamiltonian $h$ is constant along
the Hamiltonian trajectories. Moreover,
the Hamiltonian flow preserves the canonical 2-form
$\dr x\wedge\dr\xi$ and, consequently, the canonical measure
$\,\dr x\,\dr\xi\,$ on $\Bbb R^{2d}$ \cite{A,CFS}.

\subhead
{2. Hamiltonian flows on manifolds without boundary}
\endsubhead
If we identify the phase space $\Bbb R^{2d}$ with the
cotangent bundle $T^*\Bbb R^d$ (which simply means that
the choice of coordinates $x$ determines, in a standard way,
the choice of coordinates $\xi$) then the solution
$\xi^*(t;y,\eta)$ behaves under change of coordinates
as a covector over the point $\xi^*(t;y,\eta)$.
Therefore the above construction can be generalized to the
case where the Hamiltonian is a function on the cotangent
bundle $T^*M$ over a smooth $d$-dimensional manifold $M$
without boundary.

If $M$ is a manifold and
$h$ is a smooth function on $T^*M$ then the Hamiltonian equations
are understood in local coordinates, $x=(x_1,\ldots,x_d)$ and
$y=(y_1,\ldots,y_d)$ being local coordinates on (or points of) $M\,$,
$\xi=(\xi_1,\ldots,\xi_d)$ and
$\eta=(\eta_1,\ldots,\eta_d)$ being the dual
coordinates on (or points of) the fibres $T_x^*M$ and
$T_y^*M$ respectively. In this case the solution (2.1) is
a smooth trajectory in $T^*M$. As before, the Hamiltonian flow
preserves the Hamiltonian, the canonical 2-form
and canonical measure on $T^*M$ (which are defined
as $\dr x\wedge\dr\xi$ and $\,\dr x\,\dr\xi\,$
in every local coordinate system).

\subhead
{3. Homogeneous Hamiltonian flows}
\endsubhead
Let $h$ be a function on $T^*M$ satisfying the following
conditions:
\roster
\item $h(x,\xi)>0$ whenever $\xi\ne0$;
\item $h$ is infinitely smooth outside $\{\xi=0\}$;
\item $h$ is positively homogeneous of degree 1 in $\xi$,
that is, $h(x,\lambda\xi)=\lambda\,h(x,\xi)$ for all positive
$\lambda$ and all $(x,\xi)\in T^*M$.
\endroster
Clearly, such a function is not smooth at $\{\xi=0\}$.
However, we can consider the Hamiltonian
trajectories and Hamiltonian flow generated by $h$
on the invariant sets
$$
T'M\ :=\ \{\,(x,\xi)\in T^*M:h(x,\xi)\ne0\,\}\
=\ \{\,(x,\xi)\in T^*M:\xi\ne0\,\}
$$
and
$$
S^*M\ :=\ \{\,(x,\xi)\in T'M:h(x,\xi)=1\,\}\,.
$$

If conditions (1)--(3) are fulfilled then the Hamiltonian
flow has the following additional properties (see, for example,
\cite{SV2}):
\roster
\item"(i)"
the solutions $x^*$ and $\xi^*$
of the Hamiltonian system are positively homogeneous in $\eta$
of degrees 0 and 1 respectively;
\item"(ii)" $\dot x^*$ does not vanish, that is, the
trajectory $x^*(t)\subset M$ cannot stop;
\item"(iii)"
the Hamiltonian flow preserves the canonical 1-form $\xi\cdot\dr x$.
\endroster

In the theory of elliptic (pseudo)differential operators
one usually deals with the homogeneous Hamiltonian flow
generated by the principal symbol to the power $1/m\,$,
$\,m$ being the order of the operator under consideration
(see, for example \cite{MS, SV2}).

\example{Example 2.1}
Let $M$ be a Riemannian manifold
and $h(x,\xi)=|\xi|_x\,$, where $|\xi|_x$ is
the length of the covector $\xi\in T_x^*M\,$. Then
$\,x^*(t)$ are geodesics and
$$
\xi_k^*(t)\ =\ \sum_{j=1}^d
g_{jk}\bigl(x^*(t)\bigr)\,\dot x^*_j(t)\,,
$$
where $\,\{g_{jk}\}\,$ is the metric tensor and
$\,\dot x^*=\{\dot x^*_1,\ldots,\dot x^*_d\}\,$ is the
tangent vector at the point $\,x^*\,$. The homogeneous
Hamiltonian flow generated by this Hamiltonian is called
a {\sl geodesic flow\/}. Note that it is more common to
define geodesics with the use of the Hamiltonian
$$
|\xi|_x^2\ =\ \sum_{j,k=1}^d g^{jk}(x)\,\xi_j\,\xi_k\,,
$$
where $\{g^{jk}\}:=\{g_{jk}\}^{-1}$.
Since the Hamiltonian is constant along the trajectories,
this simply means that the parameter $t$ is chosen in a
different way; if $h(x,\xi)=|\xi|_x\,$ then the geodesics
are parameterized by their length.
\endexample

\head
{3. Homogeneous Hamiltonian billiards}
\endhead
Throughout this section we assume that  the Hamiltonian
satisfies the conditions (1)--(3) of subsection 2.3.
The definitions and results quoted below can be found
in \cite{SV1, SV2}.

\subhead
{1. Billiard trajectories}
\endsubhead
Let $M$ be a smooth $d$-dimensional manifold with boundary
and $h$ be a homogeneous Hamiltonian on $T^*M$.
Near $\partial M$ we shall use special coordinates
$x=(x',x_d)$ such that $x'\in\Bbb R^{d-1}\,$,
$\partial M=\{x_d=0\}$ and $x_d>0 $ for points
in the interior of $M$. Then $\xi=(\xi',\xi_d)$,
where $\xi_d$ is the so-called conormal component of $\xi$.

Let $(y,\eta)\in T'(M\setminus\partial M)$, that is,
$y\not\in\partial M\,$. Assume that, at some time $t=\tau$,
the Hamiltonian trajectory (2.1) originating from $(y,\eta)$
hits the boundary at some point
$$
(x^*(\tau-0),{\xi}^*(\tau-0))\ \in\ T'_{\partial M}M
$$
(in other words, $x^*(\tau-0)\in\partial M$). Then
$(x^*(\tau-0),{\xi}^*(\tau-0))$
is said to be a {\it point of incidence\/}. At a
point of incidence we always have
$$
h_{\xi_d}(x^*(\tau-0),\xi^*(\tau-0))\ \le\ 0\,.
$$
The incoming trajectory is said to be {\it transversal} if
$$
h_{\xi_d}(x^*(\tau-0),\xi^*(\tau-0))\ <\ 0\,.
$$

\definition{Definition 3.1}
We say that
$(x^*(\tau+0),{\xi}^*(\tau+0))\in T'_{\partial M}M$ is a
{\it point of reflection} if
\roster
\item"(i)"
$\ h(x^*(\tau-0),\xi^*(\tau-0))=h(x^*(\tau+0),\xi^*(\tau+0))\,$,
\item"(ii)"
$\ x^*(\tau-0)=x^*(\tau+0)\,$,
\item"(iii)"
the covectors $\xi^*(\tau-0)$ and $\xi^*(\tau+0)$ differ only in
their conormal  component,
\item"(iv)"
$\ h_{\xi_d}(x^*(\tau+0),\xi^*(\tau+0))>0\,$.
\endroster
Since $\,h_{\xi_d}(x^*(\tau+0),\xi^*(\tau+0))>0\,$,
the Hamiltonian trajectory originating from the point
$(x^*(\tau+0),{\xi}^*(\tau+0))$ is well defined. It is
called a {\it reflected  trajectory\/}.
\enddefinition

\remark{Remark 3.2}
Note that, generally speaking, there may be several points
of reflection and reflected trajectories corresponding
to one point of incidence.
\endremark

\definition{Definition 3.3}
The trajectory obtained by consecutive transversal
reflection is called a {\it billiard trajectory\/}.
\enddefinition

A billiard trajectory originating from
an interior point $(y,\eta)$ may not be defined for
all $t\in\Bbb R\,$ for one of the following reasons:
\roster
\item"(i)"
at some moment the incoming part of the trajectory
is not transversal,
\item"(ii)"
the trajectory experiences an infinite number of
reflections in a finite time.
\endroster
In the first case the trajectory is said to be {\it grazing\/}
and in the second case it is called a {\it dead-end\/}.

\subhead
{2. Billiard flows}
\endsubhead
If in our special coordinates
$$
\align
(x^*(\tau-0),\xi^*(\tau-0))\ &=\ ((x^*)',0,(\xi^*)',\xi_d^-)\,,\\
(x^*(\tau+0),\xi^*(\tau+0))\ &=\ ((x^*)',0,(\xi^*)',\xi_d^+)
\endalign
$$
then the reflection law can be written as
$$
\align
h\bigl((x^*)',0,(\xi^*)',\xi_d^-\bigr)\
&=\ h\bigl((x^*)',0,(\xi^*)',\xi_d^+\bigr)\,,\\
h_{\xi_d}\bigl((x^*)',0,(\xi^*)',\xi_d^+\bigr)\
&>\ 0\,.
\endalign
$$

\definition{Definition 3.4}
We say that the {\it simple reflection condition\/}
is fulfilled if, for every
$\,(x',\xi')\in T'\partial M\,$, the function
$\,h(x',0,\xi',\cdot)\,$ has the only local (and hence
global) minimum $\,\xi_d^{\hbox{st}}\,$. We say that
the {\it strong simple reflection condition\/} is fulfilled
if, in addition,
$\,h_{\xi_d\xi_d}(x',0,\xi',\xi_d^{\hbox{st}})\ne0\,$.
\enddefinition

If the simple reflection condition is fulfilled then
$$
\text{sign}\,h_{\xi_d}(x',0,\xi',\xi_d)
\ =\ \text{sign}\,(\xi_d-\xi_d^{\hbox{st}})\,,
\qquad\forall\xi_d\in\Bbb R\,.
$$
Therefore for every transversal incoming trajectory there
exists the only reflected trajectory and we have
$$
\xi_d^-\ <\ \xi_d^{\hbox{st}}\ <\ \xi_d^+\,.
$$
Under the simple reflection condition one can define
the group of shifts along the billiard trajectories,
which is called the {\it billiard flow\/}.
The billiard flow is defined for all $t\in\Bbb R\,$ on
a set of full measure \cite{CFS} (that is, the set
of starting points of grazing and dead-end trajectories is
of measure zero) and has the same properties as
the Hamiltonian flow (see subsections 2.2 and 2.3).

\example{Example 3.5}
The Hamiltonian $h(x,\xi)=|\xi|_{x}$ on a Riemannian
manifold with boundary generates so-called
{\it geodesic billiard\/}. In this case $x^*(t)$ are
geodesics (see Example 2.1) and the reflection law takes
the usual form: the angle of incidence is equal to the angle
of reflection. Clearly, the geodesic billiard satisfies
the strong simple reflection condition.
\endexample

\subhead
{3. Hamiltonian billiards with
nonnegative   Hamiltonian  curvature}
\endsubhead
If the simple reflection condition is fulfilled then
the function
$$
\bold k(x',\xi')=
\{h_{\xi_d},h\}|_{x_d=0,\,\xi_d=\xi_d^{st}(x',\xi')}
$$
is said to be the {\it Hamiltonian  curvature\/} of
$\partial M$. Here $(x',\xi')\in T'\partial M $ and
$\{\cdot,\cdot\}$ are the Poisson brackets.

\definition{Definition 3.6}
We say that the Hamiltonian billiard (or billiard flow)
is {\it convex\/} if
the strong simple reflection condition is fulfilled and
$$
\bold k(x',\xi')\ge 0\,,
\qquad \forall\,(x',\xi')\in T'\partial M\,.\tag3.1
$$
\enddefinition

\example{Example 3.7}
If $M$ is a domain in the Euclidean space and
$h(x,\xi)=|\xi|$ then (3.1) is equivalent
to the usual definition of convexity.
\endexample

In Section 6 we shall use the following result.

\proclaim{Lemma 3.8 \cite{SV2, Lemma 1.3.17}}
In a convex billiard
there are no grazing and dead-end trajectories.
\endproclaim

\head
{4. Periodic and absolutely periodic trajectories}
\endhead

\subhead
{1. Homogeneous Hamiltonian and billiard flows}
\endsubhead
Throughout this subsection we assume that
\roster
\item"(i)"
the Hamiltonian is defined on $T'M$ and
satisfies the conditions (1)--(3) of subsection 2.3;
\item"(ii)" either $\partial M=\emptyset$ or the
corresponding homogeneous Hamiltonian billiard
satisfies the simple reflection condition.
\endroster
For the
sake of convenience we shall regard homogeneous Hamiltonian
flows on manifolds without boundary as a particular
case of homogeneous billiard flows
and assume that, by definition, the Hamiltonian flows on
manifolds without boundary are convex.

\definition{Definition 4.1}
Let $\,T>0\,$.
A trajectory
$\,\bigl(x^*(t;y_0,\eta_0),\xi^*(t;y_0,\eta_0)\bigr)\,$
and its starting point
$(y_0,\eta_0)\in T'(M\setminus\partial M)$
are said to be
\roster
\item
{\it $T$-periodic\/} if
$\,\bigl(x^*(T;y_0,\eta_0),\xi^*(T;y_0,\eta_0)\bigr)
=(y_0,\eta_0)\,$;
\item
{\it absolutely $T$-periodic\/} if the function
$$
|x^*(T;y,\eta)-y|^2 +|\xi^*(T;y,\eta)-\eta |^2\tag4.1
$$
of the variables $(y,\eta)\in T'(M\setminus\partial M)$
has an infinite order zero at $(y_0,\eta_0)\,$;
\item
{\it absolutely $(T,l)$-periodic\/} if they are absolutely
$\,T$-periodic and the trajectory hits the boundary
$\,l\,$ times as $\,t\in(0,T)\,$;
\item
{\it periodic\/} if they are $T$-periodic for some $T>0\,$;
\item
{\it absolutely periodic\/} if they are absolutely $T$-periodic
for some $T>0\,$
\endroster
\enddefinition

We shall denote the sets of periodic,
absolutely periodic, absolutely $T$-periodic points and
absolutely $(T,l)$-periodic points lying in
$S^*(M\setminus\partial M)$ by
$\Pi$, $\Pi^a$, $\Pi^a_T$ and $\Pi_{T,l}^a$ respectively.

The following two lemmas suggest that
the points lying in a path-connected component of $\Pi$
have a common period. However, we do not know whether
it is true in the general case.

\proclaim{Lemma 4.2}
Let the manifold $M$ and Hamiltonian $h$
belong to a quasianalytic class $C(m_n)\,$,
and let $\Omega$ be a path-connected subset of
$\Pi^a\,$. Then either $\,\Omega\cap\Pi_{T,l}^a=\emptyset\,$
or $\Omega\subset\Pi_{T,l}^a\,$.
\endproclaim

\proclaim{Lemma 4.3}
If $\gamma$ is a smooth path in $\Pi$ and
the period $T(y,\eta)$ is a continuous function of
$(y,\eta)\in\gamma$ then $T(y,\eta)$ is constant
on $\gamma\,$.
\endproclaim

Note that in Lemma 4.2 the set $\Pi^a$ may be disconnected
even if $M$ is connected. Indeed, if there exist grazing
or dead-end trajectories then, generally speaking, the
billiard flow is well defined for $t\in[0,T]$ only on a
subset of $T'M\,$, which may well be disconnected.

The following lemma implies, in particular, that the set
$\Pi\setminus\Pi^a$ is of measure zero.

\proclaim{Lemma 4.4 \cite{SV1, SV2}}
The set of points which are $T$-periodic but
not absolutely $T$-periodic for some $T>0$ is
of measure zero.
\endproclaim

\proclaim{Theorem 4.5}
Let $M$ be a compact connected manifold
and let the billiard flow generated
by a Hamiltonian $h$ be convex. If $M$ and $h$
belong to a quasianalytic class $C(m_n)$ then
the existence of one absolutely $(T,l)$-periodic billiard
trajectory implies that all trajectories are $(T,l)$-periodic.
\endproclaim

Theorem 4.5 was proved in \cite{SV2} in the analytic case.
If $\partial M=\emptyset$ then it takes the following form.

\proclaim{Theorem 4.5$'$}
Let $M$ be a compact connected manifold without boundary.
If $M$ and $h$ belong to a quasianalytic class $C(m_n)$ then
the existence of one absolutely $T$-periodic Hamiltonian
trajectory implies that all trajectories are $T$-periodic.
\endproclaim

Theorem 4.5 implies the following important corollary.

\proclaim{Corollary 4.6}
Let $M$ be a compact connected manifold, $\partial M\ne\emptyset$,
and let the billiard flow generated by a Hamiltonian
$h$ be convex. If $M$ and $h$ belong to a quasianalytic class
$C(m_n)$  and $\bold k\not\equiv0$ then the set of periodic
points is of measure zero.
\endproclaim

The next two theorems show that the quasianalyticity
condition in Theorems 4.5 and 4.5$'$ cannot be removed.

\proclaim{Theorem 4.7}
If the class $C(m_n)$ is not quasianalytic then there exists
a Riemannian $C(m_n)$-manifold $M$ without boundary such
that the geodesic flow on $M$ satisfies the following
conditions:
\roster
\item
for some $T$
the set $\,\Pi^a_T$ has a positive measure,
\item
but there are nonperiodic trajectories; moreover,
$\,S^*M\setminus\Pi$ is a set of positive measure.
\endroster
\endproclaim

\proclaim{Theorem 4.8}
If the class $C(m_n)$ is not quasianalytic then there exists
a Riemannian $C(m_n)$-manifold $M$ with boundary such that
\roster
\item the geodesic billiard on $M$ is convex with
$\bold k\equiv0$,
\item for some $T$ and $l$ the set $\Pi_{T,l}^a$ has a
positive measure,
\item but there are nonperiodic billiard trajectories;
moreover, $S^*M\setminus\Pi$ is a set of positive measure.
\endroster
\endproclaim

One can generalize the definition of absolutely periodic
points in the following way.

\definition{Definition 4.9}
Let $T^*(y,\eta)$ be a smooth function defined in a
neighbourhood of $(y_0,\eta_0)$. We say that the point
$(y_0,\eta_0)$ is absolutely $T^*$-periodic if the function
$$
|x^*(T^*(y,\eta);y,\eta)-y|^2
+|\xi^*(T^*(y,\eta);y,\eta)-\eta |^2\tag4.2
$$
of variables $(y,\eta)$ has an infinite order zero at
$(y_0,\eta_0)$.
\enddefinition

However, in the quasianalytic case this is equivalent
to Definition 4.1. Indeed, if the manifold and functions
$\,h\,$ and $\,T^*\,$ are quasianalytic then the
function (4.2) is also quasianalytic. If it has
an infinite order zero at $(y_0,\eta_0)$ then it is
identically equal to zero in a neighbourhood of
$(y_0,\eta_0)$, which means that all points $(y,\eta)$
of this neighbourhood are periodic with period
$T^*(y,\eta)$. Now Lemma 4.3 implies that
$T^*\equiv T^*(y_0,\eta_0)$, and therefore
the point $(y_0,\eta_0)$ is absolutely
$T^*(y_0,\eta_0)$-periodic.

\subhead
{2. Branching Hamiltonian billiards}
\endsubhead
If the simple reflection condition is not fulfilled
then the corresponding billiard is said to be {\sl
branching}. In this case there may exist infinitely many
billiard trajectories originating from a fixed point
$(y,\eta)\in T'(M\setminus\partial M)$, moreover, the
set of these trajectories may well be uncountable.
Therefore Theorem 4.5, as it is stated above, is
unlikely to be true even for the simplest branching
billiards.

For homogeneous branching billiards it is possible to prove
that the set of starting points of grazing trajectories is
of measure zero \cite{SV2}. But, even in the analytic case,
the measure of the set of starting points of dead-end
trajectories may be positive \cite{SV1}.

One can classify trajectories by the type
of their reflections, introduce the notion of periodic
and absolutely periodic trajectories and prove
statements similar to Lemmas 4.2-4.4 \cite{SV2}.
However, it is of a little interest in applications
unless we have effective sufficient conditions
for the set of starting points of dead-end
trajectories to be of measure zero. To the best
of our knowledge, the only result in this
direction (without requiring simple reflection)
was obtained in \cite{V1}.

\subhead
{3. Nonhomogeneous Hamiltonians}
\endsubhead
From the geometric point of view, the nonhomogeneous Hamiltonian
flows are more difficult to study because they do not preserve
the canonical one form.
If the Hamiltonian $h$ is not homogeneous then one has to
consider the restriction of the Hamiltonian flow to
$$
\Sigma_\lambda\ :=\ \{(x,\xi)\in T^*M:h(x,\xi)=\lambda\}
$$
for each fixed $\lambda$ separately.
If $\lambda$ is not a critical value of $h$ then $\Sigma_\lambda$
is a smooth $(2d-1)$-dimensional submanifold, and one can define
the sets $\Pi\subset\Sigma_\lambda$ and $\Pi^a\subset\Sigma_\lambda$
in the same way as for homogeneous flows (with
(4.1) being considered as a function on $\Sigma_\lambda$).

In the nonhomogeneous case the problems discussed in
Introduction become much more difficult.
Not only may the answer depend on $\lambda$,
but also the structure of the sets $\Pi$ and $\Pi^a$ on
a fixed energy surface $\Sigma_\lambda$ may be more complicated.
The following simple observation shows that
Lemma 4.3 does not necessarily hold for the nonhomogeneous flows.

\example{Example 4.10} Assume that zero is not a critical value
of the Hamiltonian $h$ and let $h_g:=g\,h$, where $g$ is
a smooth strictly positive function. Then $h_g$ vanishes 
on $\Sigma_0$ and the Hamiltonian trajectories $(x_g^*,\xi_g^*)$
of $h_g$ lying on $\Sigma_0$ are defined by the equations
$$
\alignat2
\dot x_g^*(t;y,\eta)
&=g(x_g^*,\xi_g^*)\,h_\xi(x_g^*,\xi_g^*),\qquad
&&x_g^*(0;y,\eta)=y\,,\\
\dot\xi_g^*(t;y,\eta)&
=-g(x_g^*,\xi_g^*)\,h_x(x_g^*,\xi_g^*),\qquad
&&\xi_g^*(0;y,\eta)=\eta\,.
\endalignat
$$
This implies that
$x_g^*(t;y,\eta)=x^*(f(t;y,\eta);y,\eta)$
and $\xi_g^*(t;y,\eta)=\xi^*(f(t;y,\eta);y,\eta)$,
where
$$
f(t;y,\eta)=\int_0^t
g\left(x_g^*(s;y,\eta),\xi_g^*(s;y,\eta)\right)\,
\dr s\,.
$$
Therefore every point $(y,\eta)\in\Pi\subset\Sigma_0$
is periodic with respect
to the Hamiltonian flow generated by $h_g$ and its
periods $T$ and $T_g$ are related as follows
$$
T(y,\eta)\ =\ f(T_g(y,\eta);y,\eta)\,.
$$
Clearly, the period $T_g$ may vary from one point to
another even if $T$ is constant.
\endexample

In \cite{PP1, PP2} the authors proved an analogue of Lemma 4.4 and,
under certain additional restrictions, the analytic version of Theorem
4.5$'$  for a class of nonhomogeneous Hamiltonians,
including Hamiltonians of the form $h(x,\xi)=|\xi|^2+V(x)$.
The latter result is likely to remain valid in the
quasianalytic case.

\head
{5. Geodesic flows on surfaces of revolution}
\endhead

In this section we consider the geodesic flow
on a 2-dimensional surface of revolution provided with
the standard metric. In the first subsection
we write down the differential equations for geodesics in an explicit
form. We use the arguments suggested by D\. Vassiliev in
\cite{V2},
where 2-dimensional analytic manifolds whose geodesics are closed
with the same length were described (see also \cite{B, Chapter~4}).
In the second subsection we study the set
of absolutely periodic points.

\subhead
{1. Differential equations for geodesics}
\endsubhead
Let
$-\frac{\pi}{2}\le\theta\le\frac{\pi}{2}\,$,
$\,0\le\varphi\le2\pi\,$.
Consider a surface of revolution $M\subset\Bbb R^3$
defined by
$$
q_1=\cos\theta\,\sin\varphi,\quad
q_2=\cos\theta\,\cos \varphi,\quad
q_3=\int_{0}^{\theta}\sqrt{[1+f(\psi)]^2-\sin^2\psi\,}\;
\dr\psi,\tag5.1
$$
where $q_k$ are coordinates in $\Bbb R^3$ and
$f(\theta)$ is a smooth function such that
\roster
\item
$\,f(0)=0\,$,
\item
$\,1+f(\psi)>|\sin\psi|\,$ for all
$\,\psi\in(-\frac{\pi}{2},\frac{\pi}{2})\,$,
\item
$f$ can be extended to a smooth function on $\Bbb R$ which
is even with respect to the points $-\frac{\pi}{2}$ and
$\frac{\pi}{2}$ in some neighbourhoods of these points,
\item
the Taylor expansions of the function
$[1+f(\psi)]^2-\sin^2\psi$ at the points
$-\frac{\pi}{2}$ and $\frac{\pi}{2}$
start with terms of degree $2q_-$ and
$2q_+$ respectively, where $q_-$ and
$q_+$ are arbitrary odd positive integers.
\endroster

In view of (2), the surface $M$ is of the same
smoothness as the function $\,f\,$ outside the poles
$\{\theta=\pm\frac\pi2\}\,$.
The conditions (3) and (4) imply that the same is
true in a neighbourhood of the poles. Indeed, one
can easily prove that, under conditions (3) and (4),
$M$ can be defined in a neighbourhood of a pole by the
equation $z=F(r)\,$, where $r=\sqrt{x^2+y^2}$ and $F$ is
an even function which belongs to the same class $C(m_n)$
as $f\,$.

\example{Example 5.1}
If $f\in C_0^\infty(-\frac{\pi}{2},\frac{\pi}{2})$
and $f$ satisfies the conditions (1) and (2) then $f$
also satisfies (3) and (4).
\endexample

The geodesic flow on a 2-dimensional surface in $\Bbb R^3$
can be interpreted as the motion of a
particle, with velocity and mass equal to one,
in the field of inertial forces. It is
described by the Euler--Lagrange equations
$$
\frac{\dr\hfill}{\dr t}
(\partial L/\partial\dot x_i)-
\partial L/\partial x_i\ =\ 0,\qquad i=1,2,
$$
where $x_i=x_i(t)\,$, $\,\dot x_i=\dot x_i(t)$ is the tangent
vector and
$$
L=L(x_1,x_2,\dot x_1,\dot x_2)\ :=\
\sum_{k=1}^3\left(\frac{\dr\hfill}{\dr t}\,q_k(x_1,x_2)\right)^2
$$
is the Lagrangian, $x_i$ being coordinates on
the surface (see, for example \cite{DFN}).
The Euler--Lagrange equations imply that the function
$$
I_1(x_1,x_2,\dot x_1,\dot x_2)
\ :=\ \dot x_1\,(\partial L/\partial\dot x_1)
+\dot x_2\,(\partial L/\partial\dot x_2)
-L(x_1,x_2,\dot x_1,\dot x_2)
$$
is constant along the trajectories $(x(t),\dot x(t))$.
Note that this function turns into the corresponding Hamiltonian
if we replace $\dot x_i$ with $\sum_k\,g^{ik}\xi_k$, where
$\{g^{ik}\}$ is the metric tensor.

If $M$ is defined as above and
$\,x_1:=\theta\,$, $\,x_2:=\varphi\,$ then
$$
L(\theta,\dot\theta,\varphi)\ =\ I_1(\theta,\dot\theta,\varphi)
\ =\ \dot\varphi^2\cos^2\theta+\dot\theta^2\,
[1+f(\theta)]^2\,,
$$
and, by the second Euler--Lagrange equation,
the function
$$
I_2(\theta,\dot\varphi)\ =\ \dot\varphi\,\cos^2\theta
$$
is constant along the trajectories.

Having found two motion integrals $I_1$ and $I_2$ for
the Euler--Lagrange equations, we can write down equations
for the geodesics in
an explicit form. Assume that
\roster
\item "(i)"
the starting point of a geodesic
$\bigl(\varphi^*(t),\theta^*(t)\bigr)$ lies on
the equator, that is, $\,\theta^*(0)=0\,$ and so
$$
I_1\ =\ \dot\varphi^*(t)\,\cos^2\theta^*
(t)\ =\ \dot\varphi^*(0)\,;\tag5.2
$$
\item "(ii)"
the geodesic $\bigl(\varphi^*(t),\theta^*(t)\bigr)$ is
parametrized by its length, that is,
$$
I_2\ =\ (\dot\varphi^*(t))^2\cos^2\theta^*(t)
+(\dot\theta^*(t))^2\,
\bigl[1+f(\theta^*(t))\bigr]^2\ =\ 1\,.\tag5.3
$$
\endroster
Clearly, $\dot\varphi^*(0)=\cos\alpha\,$, where
$\alpha\in[0,\frac\pi2]$ is the angle
between the geodesic and equator at the starting point.
Therefore we can rewrite (5.2) as
$$
\dot\varphi^*(t)\ =\
\frac{\cos\alpha\hfill}{\cos^2\theta^*(t)}\,. \tag5.4
$$
Now (5.3) and (5.4) imply
$$
|\dot\theta^*(t)|\ =\
\frac{\sqrt{\cos^2\theta^*(t)-\cos^2\alpha}}
{\bigl[1+f(\theta^*(t))\bigr]\cos\theta^*(t)}\,.\tag5.5
$$
Note that $|\theta^*(t)|\le\alpha$ and therefore
$\cos^2\theta^*(t)\ge\cos^2\alpha\,$. Indeed,
if $\dot\theta^*(t)=0\,$ then, in view of (5.3) and (5.4),
$\,|\theta^*(t)|=\alpha\,$. In particular, if
$\theta^*(t)$ has a local maximum at $t_+$ and
a local minimum at $t_-$ then
$$
\theta^*(t_+)\,=\,-\theta^*(t_-)\,=\,\alpha\,.\tag5.6
$$

The last observation together with the fact that
$\,(\dot\theta^*(t))^2=1-\cos^2\alpha\,$ whenever
$\theta^*(t)=0$ imply that the function $\theta^*(t)$
is periodic. Moreover, if
$$
\dot\theta^*(t_1)=\dot\theta^*(t_2)=\dot\theta^*(t_3)=0\,,
\qquad t_1<t_2<t_3\,,
$$
and $\,\dot\theta^*(t)\ne0\,$ for
$\,t\in(t_1,t_2)\cup(t_2,t_3)\,$
then its period coincides with $\,t_3-t_1\,$.
In view of (5.5) and (5.6) we have
$$
t_2-t_1\ =\ t_3-t_2\
=\ \int_{-\alpha}^\alpha\frac{[1+f(\theta)]\,\cos\theta}
{\sqrt{\cos^2\theta-\cos^2\alpha}}\;\dr\theta\
=\ \pi+\int_{-\alpha}^\alpha\frac{\tilde f(\theta)\,\cos\theta}
{\sqrt{\cos^2\theta-\cos^2\alpha}}\;\dr\theta\,,
$$
where
$$
\tilde f(\theta)\ :=\ \frac12\,(f(\theta)+f(-\theta))\,.
$$
Therefore
$$
\theta^*(t+T_*(\alpha))\ =\ \theta^*(t)\tag5.7
$$
for all $t\in\Bbb R\,$, where
$$
T_*(\alpha)\ :=\
2\pi+2\int_{-\alpha}^\alpha\frac{\tilde f(\theta)\,\cos\theta}
{\sqrt{\cos^2\theta-\cos^2\alpha}}\;\dr\theta
$$
is the period of $\theta^*(t)\,$.

Similarly, in view of (5.4)--(5.6),
$$
\varphi^*(t_2)-\varphi^*(t_1)
\ =\ \varphi^*(t_3)-\varphi^*(t_2)\
=\ \cos\alpha\int_{-\alpha}^\alpha\frac{1+f(\theta)}
{\cos\theta\,\sqrt{\cos^2\theta-\cos^2\alpha}}\;\dr\theta\,.
$$
Using
$$
\frac{1}{\cos\theta\,\sqrt{\cos^2\theta-\cos^2\alpha}}
\ =\ \frac{\dr\hfill}{\dr\theta}\left(\arctan\left(
\frac{\cos\alpha\sin\theta}
{\sqrt{\cos^2\theta-\cos^2\alpha}}\right)\right)
$$
we obtain
$$
\varphi^*(t_3)-\varphi^*(t_1)\ =\ 2\pi+R(\alpha)\,,\tag5.8
$$
where
$$
R(\alpha)\ :=\
2\cos\alpha\int_{-\alpha}^\alpha\frac{\tilde f(\theta)}
{\cos\theta\,\sqrt{\cos^2\theta-\cos^2\alpha}}\;\dr\theta\,.
$$
Since the function $\theta^*(t)$ is $T_*(\alpha)$-periodic,
(5.4) and (5.8) imply that for all $t\in\Bbb R\,$
$$
\varphi^*(t+T_*(\alpha))
\ =\ \varphi^*(t)+R(\alpha)\pmod{2\pi}\,.\tag5.9
$$

\subhead
{2. Absolutely periodic points}
\endsubhead
In the previous subsection we have proved that any geodesic
starting at the equator can be written as
$(\varphi^*(t),\theta^*(t))\,$, where the functions
$\varphi^*(t)$ and $\theta^*(t)$ satisfy the differential
equations (5.4), (5.5). If we identify vectors and covectors
as in Example 2.1 then
the corresponding Hamiltonian
trajectory has the form
$(\varphi^*(t),\theta^*(t),\dot\varphi^*(t),\dot\theta^*(t))$.
In particular, the equator is a $2\pi$-periodic geodesic
(with $\dot\varphi^*(0)=1$), and the corresponding Hamiltonian
trajectory is $\,\Gamma(t)=(t+\varphi^*(0),0,1,0)\,$,
$t\in\Bbb R\,$. The following lemma gives a necessary
and sufficient condition for this trajectory to be
absolutely $2\pi$-periodic.

\proclaim{Lemma 5.2}
The Hamiltonian trajectory $\Gamma(t)$ is
absolutely $2\pi$-periodic if and only if all even derivatives
of the function  $f(\theta)$ vanish at $\theta=0\,$, that is,
$$
f^{(2k)}(0)=0\,,\qquad k=0,1,2,\ldots\tag5.10
$$
\endproclaim

\demo{Proof}
Let
$\,
\left(\theta^*(t;\theta,\varphi,\dot\theta,\dot\varphi),
\varphi^*(t;\theta,\varphi,\dot\theta,\dot\varphi)\right)
\,$
be a geodesic starting at the point
$(\varphi,\theta,\dot\varphi,\dot\theta)\,$
and $\alpha^*(\theta,\dot\varphi)$ be the angle of
intersection of this geodesic with equator (here we
consider $\,\varphi,\theta,\dot\varphi,\dot\theta\,$
as independent variables).
Assume, for the sake of definiteness, that
$\varphi=0\,$. We have to prove that
the functions
$$
\varphi^*(2\pi;\varphi,\theta,\dot\varphi,\dot\theta)
-\varphi\,,\
\theta^*(2\pi;\varphi,\theta,\dot\varphi,\dot\theta)
-\theta\,,\
\dot\varphi^*(2\pi;\varphi,\theta,\dot\varphi,\dot\theta)
-\dot\varphi\,,\
\dot\theta^*(2\pi;\varphi,\theta,\dot\varphi,\dot\theta)
-\dot\theta
$$
have infinite order zeros at the point $\{0,0,1,0\}\,$
if and only if (5.10) holds true.

According to (5.4)
$$
\cos\alpha^*=\dot\varphi\,\cos^2\theta\tag5.11
$$
and, by (5.7),
$$
\theta\ \equiv\ \theta^*(T_*(\alpha^*);
\varphi,\theta,\dot\varphi,\dot\theta)\,.\tag5.12
$$
If the function
$\,
\theta^*(2\pi;\varphi,\theta,\dot\varphi,\dot\theta)-\theta
\,$
has an infinite order zero at $\{0,0,1,0\}\,$ then,
in view of (5.5) and (5.12), the function
$\,T_*(\alpha^*)-2\pi\,$ has an infinite order zero
at $\{\theta=0,\dot\varphi=1\}\,$. This fact and
(5.11) imply that $T_*(\alpha)-2\pi$ has an infinite order
zero at $\alpha=0\,$, which is equivalent to (5.10).

Assume now that all even derivatives of $f$
vanish at $\theta=0\,$. Then the functions $R(\alpha)$
and $\,T_*(\alpha^*)-2\pi\,$ have infinite order zeros
at $\{\theta=0,\dot\varphi=1\}\,$, so (5.12) and
(5.9) respectively imply that
$\,
\theta^*(2\pi;\varphi,\theta,\dot\varphi,\dot\theta)-\theta
\,$
and
$\,
\varphi^*(2\pi;\varphi,\theta,\dot\varphi,\dot\theta)-\varphi
\,$
have infinite order zeros at $\{0,0,1,0\}\,$. By (5.4)
and (5.5) the same is true for the functions
$\,\dot\varphi(2\pi;\varphi,\theta,\dot\varphi,\dot\theta)
-\dot\varphi\,$ and
$\,\dot\theta(2\pi;\varphi,\theta,\dot\varphi,\dot\theta)
-\dot\theta\,$.
\qed\enddemo

\proclaim{Lemma 5.3}
If $f$ is odd on an interval $[-\alpha_0,\alpha_0]\,$,
$0<\alpha_0\le\frac{\pi}{2}\,$ then
all the geodesics intersecting the equator at an
angle $\alpha\le\alpha_0$ are $2\pi$-periodic.
In particular, if $f$ is odd on $[-\pi/2,\pi/2]$
then all geodesics are $2\pi$-periodic.
\endproclaim

\demo{Proof}
Under conditions of the lemma
$$
T_*(\alpha)\ =\ 2\pi\,,\quad
R(\alpha)=0\,,\qquad
\forall\alpha\in[-\alpha_0,\alpha_0]\,.
$$
Therefore the required result immediately follows
from (5.7) and (5.9).
\qed\enddemo

\head
{6. Proofs}
\endhead

\subhead
{1. The period function}
\endsubhead
Let $(y_0,\eta_0)$ be a $T$-periodic point of the homogeneous
Hamiltonian or billiard flow. We are going to define in a
neighbourhood of this point
a positive function $t^*(y,\eta)$ with the following properties:
\roster
\item
$t^*(y_0,\eta_0)=T$ and
$t^*$ is positively homogeneous in $\eta$ of degree 0;
\item
$t^*$ belongs to the same Carleman class $C(m_n)$ as the pair
$(h,M)$;
\item
$\nabla_{y,\eta}t^*(y,\eta)=0$ if and only if the point $(y,\eta)$
is periodic with period $t^*(y,\eta)\,$,
\item
$\nabla_{y,\eta}t^*$ has an infinite order zero at $(y,\eta)$
if and only if the point $(y,\eta)$ is absolutely
$t^*(y,\eta)$-periodic.
\endroster
Note that the function $t^*$ is not uniquely defined by (1)--(4).
Indeed, outside the set $\Pi$ we only assume (1) and (2).
Moreover, if $(y,\eta)\in\Pi$ or $(y,\eta)\in\Pi^a$ but the corresponding
period does not coincide with $t^*(y,\eta)$ then we do not impose
any restrictions on the derivatives of $\,t^*$ at $(y,\eta)\,$.

\proclaim{Lemma 6.1}
Let the manifold $M$ and Hamiltonian $h$ belong to a class
$C(m_n)$ and let $(y_0,\eta_0)$ be a $T$-periodic point.
Then there exists a positive function $t^*$ satisfying the
conditions {\rm (1)--(4)} in a neighbourhood of
$(y_0,\eta_0)\,$.
\endproclaim

The proof is based on the following technical
lemma (see \cite{SV2, Lemma~2.3.2} or \cite{T}).

\proclaim{Lemma 6.2}
Under conditions of Lemma 6.1
there exists a coordinate system in
a neighbourhood of $y_0$ such that
$\det\xi_\eta(T;y_0,\eta_0)\ne0\,$.
\endproclaim

Note that in \cite{SV2, T} the authors only proved the existence
of $C^\infty$-coordinates for which the matrix
$\det\xi_\eta(T;y_0,\eta_0)$ is nondegenerate.
However, it is clear
from the proof that in the case of $C(m_n)$-manifold one can
find $C(m_n)$-coordinates with the same property.

\demo{Proof of Lemma 6.1}
Let us choose coordinates as in Lemma 6.2 and consider
the function
$$
\Phi(t;x,y,\eta)\ =\ (x-x^*(t;y,\eta))\cdot\xi^*(t;y,\eta)
$$
defined in a neighbourhood of the point $(T;y_0,\eta_0)\,$.
Using the fact that the Hamiltonian and billiard flows preserve
the Hamiltonian and canonical  1-form on $T'M$, one can prove
that in a neighbourhood of $(T;y_0,\eta_0)$
$$
\Phi_\eta(t;x,y,\eta)=0\quad\Longleftrightarrow
\quad x=x^*(t;y,\eta)\,,\tag6.1
$$
$$
\align
\left.\Phi_x(t;x,y,\eta)\right|_{x=x^*(t;y,\eta)}
\ &=\ \xi^*((t;y,\eta)\,,\tag6.2\\
\left.\Phi_y(t;x,y,\eta)\right|_{x=x^*(t;y,\eta)}
\ &=\ -\eta\,,\tag6.3\\
\left.\Phi_t(t;x,y,\eta)\right|_{x=x^*(t;y,\eta)}
\ &=\ -h(y,\eta)\tag6.4
\endalign
$$
(see \cite{SV2, Sections 2.3, 2.4} for details).

Define
$$
\Psi(t;y,\eta)\ :=\ \Phi(t;y,y,\eta)\,.
$$
In view of (6.4),
in a neighbourhood of $(y,\eta)$ the equation
$$
\Psi(t;y,\eta)\ =\ 0
$$
has the only $t$-solution $t^*(y,\eta)\,$. Since the point
$(y_0,\eta_0)$ is $T$-periodic we have
$\Psi(t;y,\eta)=0\,$,
and therefore $t^*(y_0,\eta_0)=T\,$. Clearly, the function
$t^*$ is positively homogeneous in $\eta$ and, by Lemma 1.5,
it belongs to the class $C(m_n)\,$.

From (6.1)--(6.3) it easily follows that
$\,\nabla_{y,\eta}\Psi(t;y,\eta)=0\,$ if and only if
the point $(y,\eta)$ is $t$-periodic, and
$\,\nabla_{y,\eta}\Psi(t;y,\eta)\,$ has an infinite order
zero at $(y_0,\eta_0)$ (as a function of $(y,\eta)$) if and only
if the point  $(y_0,\eta_0)$ is absolutely $t$-periodic
(see \cite{SV2, Section~4.1} for details). Differentiating
the identity
$$
\Psi(t^*(y,\eta);y,\eta)\ =\ 0
$$
in $y$ and $\eta$ and taking into account (6.4), we see that
$\left.\nabla_{y,\eta}\Psi(t;y,\eta)\right|_{t=t^*(y,\eta)}=0\,$
if and only if $\,\nabla t^*(y,\eta)=0\,$, and
$\left.\partial_y^\alpha\partial_\eta^\beta
\Psi(t;y_0,\eta_0)\right|_{t=t^*(y_0,\eta_0)}=0\,$
for all nonzero $\alpha,\beta$ if and only if $\nabla t^*(y,\eta)$
has an infinite order zero at $(y_0,\eta_0)$. This implies (3) and (4).
\qed\enddemo

\subhead
2. Proof of Lemma 4.2
\endsubhead
Let $(y_0,\eta_0)\in\Pi_{T,l}^a\,$.
Since we consider only transversal reflections,
the billiard trajectories
$\bigl(x^*(t;y,\eta),\xi^*(t;y,\eta)\bigr)$
are well defined for all $\,t\in[0,T]\,$ and
all $(y,\eta)\,$ lying in a neighbourhood of
$(y_0,\eta_0)$. Moreover, every trajectory
starting in this neighborhood hits the boundary
$\,l\,$ times as $t\in[0,T]\,$.

Let $t^*$ be the period function defined on a smaller
neighbourhood of $(y_0,\eta_0)\,$.
In view of (4) the gradient of this function has an infinite
order zero at $(y_0,\eta_0)\,$. Under conditions of
the lemma $t^*$ is quasianalytic, and therefore it is
identically equal to $T$ in this neighbourhood.
Thus, every point $(y,\eta)\in\Pi_{T,l}^a\,$
has a neighbourhood $U_{(y,\eta)}\subset\Pi_{T,l}^a\,$.

Let $\gamma$ be a path in $\Omega\,$ and,
for $\,(y,\eta)\in\gamma\,$, let
$\,T(y,\eta):=\min\{T:(y,\eta)\in\Pi_T^a\}\,$.
If $(y^{(n)},\eta^{(n)})\to(y,\eta)$ and
$T(y^{(n)},\eta^{(n)})\to T\,$ as $n\to\infty$ then,
obviously, $(y,\eta)\in\Pi_T^a\,$.
This implies that the function $T(y,\eta)$ takes its
minimal value
on $\gamma$ at some point $(y_0,\eta_0)\in\gamma$.
Consider the set $\gamma\cap\Pi_{T_0,l_0}^a\,$.
By the above, this set is open in $\gamma\,$.
On the other hand, since the restriction of the
billiard flow
$$
(y,\eta)\ \to\
\bigl(x^*(T_0;y,\eta),\xi^*(T_0;y,\eta)\bigr)
$$
to a neighbourhood of $\gamma$ is a smooth map, the set
$\gamma\cap\Pi_{T_0,l_0}^a\,$ is closed in $\gamma\,$.
Therefore $\gamma\subset\Pi_{T_0,l_0}^a\,$.
Since any two points of $\Omega$ can be joined by a
path, this implies that $\Omega\subset\Pi_{T,l}$
provided that $\Omega\cap\Pi_{T,l}\ne\emptyset\,$.

\subhead
3. Proof of Lemma 4.3
\endsubhead
Under conditions of the lemma, in a neighbourhood
of every point $(y_0,\eta_0)$ the function $T(y,\eta)$
coincides with the restriction of the period function
$t^*$ to $\gamma\,$. Now (3) implies that
$T(y,\eta)$ is locally and hence globally constant on
$\gamma\,$.

\subhead
4. Proof of Lemma 4.4
\endsubhead
The lemma immediately follows from the fact that the set of
zeros  of $\nabla_{y,\eta}$ is of the same measure
as the set of infinite order zeros.

\subhead
5. Proof of Theorem 4.5
\endsubhead
By Lemma 3.8, under conditions of the theorem
the trajectories $(x^*(t; y,\eta),\xi^*(t; y,\eta))$
are well defined for all  $\,t\in\Bbb R\,$ and
all $(y,\eta)\in T'(M\setminus\partial M)\,$.
Let $(y_0,\eta_0)\in\Pi_{T,l}^a\,$ and
$\gamma$ be a path in $T'(M\setminus\partial M)\,$
from $(y_0,\eta_0)$ to another point $(y,\eta)$.
In the same way as in the proof of Lemma 4.2
one can show that $\gamma\subset\Pi_{T,l}^a\,$.
Since the set $T'(M\setminus\partial M)\,$ is
path-connected, this implies that
$\Pi_{T,l}^a=T'(M\setminus\partial M)\,$.

\subhead
6. Proof of Corollary 4.6
\endsubhead
In view of Theorem 4.5 it is sufficient to show
that, given $l$ and $T$ we can find at least one
trajectory which is not $(T,l)$-periodic. If
$\bold k(x',\xi')\ne0$ then, by choosing the starting
point close to $(x',0,\xi',0)$, we can construct a
trajectory which experiences arbitrary many reflections
in any given time (see proof of Lemma 1.3.34 in
\cite{SV2}). This proves the corollary.

\subhead
7. Proof of Theorems 4.7 and 4.8
\endsubhead
Every nonquasianalytic class $C(m_n)$ contains
nonnegative $C_0^\infty$-functions \cite{M}. Let
$f\in C(m_n)\cap C_0^\infty(\Bbb R)\,$, $f\ge0$ and $f=0$
outside the interval $(\frac12,1)\,$. Consider the surface
of revolution $M$ defined by (5.1). Then, according to Lemma
5.3, all the geodesics intersecting the equator at an
angle $\alpha\le\frac12$ are $2\pi$-periodic, which implies
that the set $\Pi_{2\pi}^a$ is of positive measure.

On the other hand, in view of (5.7) and (5.9), a geodesics
intersecting the equator at an angle $\alpha$ is periodic
only if
$$
\multline
k\,R(\alpha)\
=\ 2k\,\cos\alpha\int_{-\alpha}^\alpha\frac{\tilde f(\theta)}
{\cos\theta\,\sqrt{\cos^2\theta-\cos^2\alpha}}\;\dr\theta\\
=\ 2k\,\cos\alpha\int_0^\alpha\frac{f(\theta)}
{\cos\theta\,\sqrt{\cos^2\theta-\cos^2\alpha}}\;\dr\theta\
=\ 0\pmod{2\pi}
\endmultline\tag6.6
$$
for some integer $k\,$.
Since $f$ has an infinite order zero at $\alpha=\frac12\,$
and $\cos\alpha$ is a decreasing function on $(0,\frac\pi2)\,$,
the function $R(\alpha)$ is strictly decreasing on an interval
$[\frac12,\alpha']\,$, $\alpha'\in(\frac12,\frac\pi2]\,$.
Therefore, for each $k\,$, (6.6) can only be true for a finite
number of points $\alpha\in[\frac12,\alpha']\,$. This implies
that the set of periodic points corresponding to the
trajectories with $\alpha\in[\frac12,\alpha']\,$ is of measure
zero. The measure of the set of starting points of
all Hamiltonian trajectories with
$\alpha\in[\frac12,\alpha']\,$ is not zero,
so $S^*M\setminus\Pi$ is a set of positive measure.

In order to prove Theorem 4.8 it is sufficient to notice that
our surface is symmetric with respect to the plane $\{x=0\}\,$.
Therefore the billiard trajectories on the manifold with
boundary $M_+=\{(x,y,z)\in M:x\ge0\}$ behave in the same way
as the geodesics on $M\,$; namely, the reflected trajectories
of the billiard flow are obtained by the reflection with
respect to the plane $\{x=0\}\,$ of the parts of geodesics
lying in the second half of $M\,$.
This implies that, for the function $f$ described above,
$\Pi_{2\pi,2}^a$ and $S^*M\setminus\Pi$ are sets of positive
measure.

The fact that $\bold k\equiv0$ is easily checked by a direct
calculation.

\enddocument